\documentclass{amsart}
\RequirePackage{float}%图片浮动
\RequirePackage{graphicx}% 插图包
\RequirePackage{amsmath, amsthm, amssymb, amsfonts,slashed}%数学环境
\RequirePackage[colorlinks=true]{hyperref}%交叉引用
\hypersetup{urlcolor=blue, citecolor=blue}%链接颜色
\RequirePackage[nameinlink]{cleveref}%智能引用
\RequirePackage{xifthen}%判断逻辑
\RequirePackage{subcaption}
\RequirePackage{indentfirst}
\RequirePackage{microtype}%文字间距优化(和中文一起)
\newenvironment{Eq}[1][]{\begin{equation}\ifthenelse{\equal{#1}{}}{}{\tag{#1}}\begin{aligned}}{\end{aligned}\end{equation}\ignorespacesafterend}
\newenvironment{Eq*}[1][]{\begin{equation*}\ifthenelse{\equal{#1}{}}{}{\tag{#1}}\begin{aligned}}{\end{aligned}\end{equation*}\ignorespacesafterend}
\numberwithin{equation}{section}
\newtheorem{theorem}{Theorem}[section]
\newtheorem{proposition}[theorem]{Proposition}

\theoremstyle{remark}%
\newtheorem{remark}{Remark}[section] 
\theoremstyle{definition}
\DeclareMathOperator{\Id}{Id}
\DeclareMathOperator{\fd}{d}	\renewcommand{\d}{\fd}
\DeclareMathOperator{\diag}{diag}

\newcommand{\Roma}[1]{\uppercase\expandafter{\romannumeral#1}}
\newcommand{\varE}{{\mathfrak E}}

\newcommand{\rH}{{\mathcal H}}

\newcommand{\rD}{{\mathcal D}}
\newcommand{\rK}{{\mathcal K}}

\newcommand{\rR}{{\mathbb{R}}}

\newcommand{\spartial}{{\slashed \partial}}

\newcommand{\kl}[1]{\mathopen{}\left#1}
\newcommand{\kr}[1]{\right#1}

\crefname{equation}{}{}%配套于cleverref
\newcommand{\Tm}[1]{Theorem \ref{#1}}

\newcommand{\Pn}[1]{Proposition \ref{#1}}

\newcounter{part0}[subsection]
\renewcommand{\part}[1][]{\noindent\refstepcounter{part0}{\bfseries Part \number\value{part0}:} #1.\par}

\newcounter{part1}[part0]

\newcounter{part2}[part1]

\title[The defocusing wave equations in two dimensional exterior region]
{The global behaviors for defocusing wave equations in two dimensional exterior region}

\author{Wei Dai}
\address{Beijing International Center for Mathematical Research\\ Peking University\\ Beijing, China}
\email{daiw21@pku.edu.cn}
\date{\today}
\begin{document}

\bibliographystyle{plain}

\begin{abstract}
We study the defocusing semilinear wave equation in $\rR\times\rR^2\backslash\rK$ with the Dirichlet boundary condition, where $\rK$ is a star-shaped obstacle with smooth boundary. 
We first show that the potential energy of the solution will decay appropriately.
Based on it, we show that the solution also pointwisely decays to $0$.
Finally, we show that the solution scatters both in energy space and the critical Sobolev space.
In general, we show that most of the conclusions obtained in \cite{MR4395159},
which hold on $\rR^{1+2}$,
remain valid on $\rR\times\rR^{2}\backslash\rK$.
\end{abstract}

\keywords{asymptotic behavior, defocusing wave equation, exterior region}

\subjclass[2020]{35L05, 35L15, 35L71}

\maketitle
%\tableofcontents

\section{Introduction}\label{Sn:I}
In this paper, we study the global asymptotic behaviors for solutions of the following defocusing semilinear wave equation with Dirichlet boundary condition:
\begin{Eq}\label{Eq:oe}%original equation
\begin{cases}
\square\phi:=\partial_t^2\phi-\Delta\phi=-|\phi|^{p-1}\phi,\qquad (t,x)\in \rR\times\rR^2\backslash\rK,\\
\phi|_{\rR\times\partial\rK}=0,\\
(\phi,\partial_t\phi)_{t=0}=(\phi_0,\phi_1),
\end{cases}
\end{Eq}
with $p>1$ and $\rK$ is a star-shaped obstacle with smooth boundary.
We shall show that most of the asymptotic behaviors which hold for the solution to the defocusing wave equation on $\rR^{1+2}$ remain valid on $\rR\times\rR^{2}\backslash\rK$.

Such a nonlinear model has drawn extensive attention in the past decades.
From the structure of \eqref{Eq:oe}, we can easily deduce that there is a conserved energy 
\begin{Eq*}
\int_{\rR^2\backslash\rK}\frac{1}{2}|\partial u|^2+\frac{1}{p+1}|u|^{p+1}\d x.
\end{Eq*}
Therefore, for general time-space with space dimension $d$, 
\eqref{Eq:oe} split naturally into two classes depending on whether or not the kinetic energy dominates the potential energy, 
where the critical exponent is given by $p=(d+2)/(d-2)$. 
Specially, for $d=2$, we see the equation always belongs to the sub-critical cases.

On $\rR^{1+2}$, the global existence of solution in energy space is well studied, 
see for example Ginibre and Velo \cite{MR786279,MR984146}.
Meanwhile, since there exists a positive definite and conserved energy,
it is expected that the solution should decay in some sense.
In early work, Glassey and Pecher \cite{MR667924}
used the energy method with respect to the conformal Killing vector field $K^0$
and showed that the time decay rate of the solution is $-1/2$ while $p>5$, 
which is consistent with the case of the homogeneous wave equation. 
They also showed a slower decay rate of solution for the range of $1+2\sqrt{2}<p\leq 5$.
From the process in that paper, it is easy to see that the vector field $K^0$ can only play a limited role for sub-conformal cases $p<5$.
To better investigate the cases for small $p$, 
Wei and Yang \cite{MR4395159} introduced classes of new vector fields 
and finally showed that the solution also enjoys the sharp decay rate $-1/2$ for $11/3<p<5$,
and still decays for $1<p\leq 11/3$.

Based on the above observations, 
it is also natural to expect that the solution should be scattered in some senses for $p$ that is not too small.
In \cite{MR667924}, the authors showed that the solution scatters in energy space when $p>p_0$ with $p_0\approx 4.15$.
And in \cite{MR4395159}, such region of $p$ has been improved to $p>2\sqrt{5}-1$.
Meanwhile, it also showed that the solution scatters in critical Sobolev space when $p>1+2\sqrt{2}$.
Besides, Ginibre and Velo \cite{MR921307} established the complete scattering
theory with data in conformal energy space for super-conformal cases $p>5$.
For a more comprehensive discussion of scattering and global dynamics for the defocusing semilinear wave equation in general dimension,
we refer interested readers to see Hidano \cite{MR1875149} and Yang \cite{MR4411879}.

Backing to the problem with Dirichlet boundary condition, we expect similar results as above. 
When the space dimension is $d=3$, 
Smith and Sogge \cite{MR1308407} first showed the global existence of smooth solution in energy critical case $p=5$ while $\rK$ is convex. 
Later, Burq, Lebeau and Planchon \cite{MR2393429} established the global existence result in energy space for equation on any bounded domain, 
which can also be extended to a general region since wave operator $\square$ has a finite propagation speed.
On the other hand, Blair, Smith and Sogge \cite{MR2566711} proved the solution is scattered while $\rK$ is star-shaped.
These authors considered the most difficult energy critical case,
yet such a process can be easily modified to the sub-critical cases at least for $p\geq 3$.

However, in the case that space dimension is $d=2$, 
even though \eqref{Eq:oe} always being sub-critical,
less is known about its behaviors.
On the one hand, 
as that in $\rR^{1+2}$, 
it is difficult to find a suitable vector field when $p$ is small.
Particularly, the presence of the obstacle $\rK$ gives additional restrictions to the vector field
and also inhibits us when obtaining pointwise estimates.
On the other hand, for example,
those Strichartz estimates used in \cite{MR4395159}
are hard to generalize in exterior region.
As discussed by Vainberg \cite{MR0415085}, Burq \cite{MR2001179} and Ralston \cite{MR552527},
the local energy estimate that required to prove the global Strichartz estimate
holds only for the Dirichlet boundary condition but not for the Neumann boundary condition.
Meanwhile, as shown in \cite{MR2584618} by Hidano, Metcalfe, Smith, Sogge and Zhou,
\cite{MR1789924} by Smith and Sogge, 
and \cite{MR2084401} by Metcalfe,
the Sobolev index $s$ can only be $1/2$ while deducing the global Strichartz estimate, 
which is not enough to show scattering.
One approach to bypass this limitation is to establish a weak version of Strichartz estimate,
such as that obtained in Smith, Sogge and Wang \cite{MR2888248}.
However, unluckily, their results can not be used in our problem.

To overcome these difficulties, we make a series of improvements.
Firstly, inspired by the vector fields used in \cite{MR4395159}, 
we construct a spherical symmetric vector field that fits our boundary condition and the shape of $\rK$ for any $p>1$.
Then, we point out that the constant that deducing pointwise estimate by embedding is uniform,
though the obstacle $\rR\times\rK$ does not scale.
Finally, instead of using the Strichartz estimate, 
we use the mixed time-space estimate directly induced by the energy inequality.
The obtained decay estimates make it possible to employ the scattering result.

To state our theorems, we make some definitions here.
Without loss of generality, we assume $\rK$ subset to the ball $B_{R}:=B(0,R)$ with radius $R\gg1$.
For writing convenience, we define $s_p:=(p-3)/(p-1)$ to be the critical Sobolev index.
We denote $a\lesssim b$ and $b\gtrsim a$ if $a\leq Cb$ for some constant $C>0$ with obvious dependence,
which may change from line to line.
We define the weighted energy norm of the initial data
\begin{Eq*}
\varE_{k,\gamma}=\sum_{l\leq k}\int_{\rR^2\backslash\rK} (1+|x|)^{\gamma+2l}(|\partial_x^{l+1}\phi_0|^2+|\partial_x^{l}\phi_1|^2)+(1+|x|)^\gamma|\phi_0|^{p+1}\d x.
\end{Eq*}
We will only use the conformal energy $\varE_{0,2}$ and the first order standard energy $\varE_{1,0}$ in this paper.

Our first result is to show the decay of the potential energy.
\begin{theorem}\label{Tm:ped}%potential energy decay
Assume $\phi$ is the solution of \eqref{Eq:oe}. 
When $1<p<5$, we have the potential energy estimate
\begin{Eq}\label{Eq:ped1}%potential energy decay 1
\int_{\rR^2\backslash\rK}(1+|t|+|x|)^{\frac{p-1}{2}}|\phi|^{p+1}\d x\lesssim \varE_{0,2}.
\end{Eq}
When $p\geq 5$, we have the potential energy estimate
\begin{Eq}\label{Eq:ped2}%potential energy decay 2
\int_{\rR^2\backslash\rK}(1+|t|+|x|)^2|\phi|^{p+1}\d x\lesssim \varE_{0,2}.
\end{Eq}
\end{theorem}

Based on the decay of the potential energy, we obtain the pointwise decay estimates of the solution.
\begin{theorem}\label{Tm:pd}%pointwise decay
For any $\varepsilon>0$, the solution $\phi$ of \eqref{Eq:oe} satisfies the pointwise decay estimate
\begin{Eq}
|\phi(t,x)|\lesssim 
\begin{cases}
(1+|t|+|x|)^{-\frac{p-1}{8}+\varepsilon}(1+||t|-|x||)^{-\frac{p-1}{8}},&p<5,\\
(1+|t|+|x|)^{-\frac{1}{2}+\varepsilon}(1+||t|-|x||)^{-\frac{1}{2}},&p\geq5,
\end{cases}
\end{Eq}
where the constant of  inequality depend on $\varepsilon$, $\varE_{0,2}$ and $\varE_{1,0}$.
\end{theorem}
\begin{remark}
Due to the lack of fundamental solution to Dirichlet-wave equation,
we can not obtain the same decay rate as that presented in \cite{MR4395159} in this article.
However, we believe that those decay rates still hold for exterior cases.
\end{remark}

Mixing the above two results, we are able to show the scattering result of \eqref{Eq:oe}.
Define $L(t)$ to be the linear wave propagation operator as follows:
\begin{Eq*}
L(t)(\phi_0,\phi_1)=(\psi,\partial_t\psi)\quad\text{satisfying}\quad L(0)=\Id,\quad\square\psi=0.
\end{Eq*}
Then, we have the result.
\begin{theorem}\label{Tm:sr}%scattering result
Assume $\varE_{1,0}$ and $\varE_{0,2}$ are finite and $p>2\sqrt{5}-1$.
Then, $\phi$ scatters in energy space, that is, there exists $(\phi_0^{\pm},\phi_1^{\pm})$ such that
\begin{Eq}\label{Eq:sr1}%scattering result
\lim_{t\rightarrow\pm\infty}\|(\phi,\partial_t\phi)-L( t)(\phi_0^{\pm},\phi_1^{\pm})\|_{\dot H^1\times L^2}=0.
\end{Eq}
Moreover, assume that $(\phi_0,\phi_1)\in \dot H^{s_p}\times \dot H^{s_p-1}$ with $p>1+2\sqrt{2}$.
Then, $\phi$ scatters in $\dot H^s\times \dot H^{s-1}$ for all $s_p\leq s< 1$, that is
\begin{Eq}\label{Eq:sr2}%scattering result
\lim_{t\rightarrow\pm\infty}\|(\phi,\partial_t\phi)-L( t)(\phi_0^{\pm},\phi_1^{\pm})\|_{\dot H^s\times\dot H^{s-1}}=0.
\end{Eq}
\end{theorem}
\begin{remark}
The homogeneous Sobolev space $\dot H^s$ on $\rR^2\backslash\rK$ is defined with the norm
\begin{Eq*}
\|f\|_{\dot H^s}=\|\Lambda^s f\|_{L_x^2(\rR^2\backslash\rK)},
\end{Eq*}
with $\Lambda=\sqrt{-\Delta_D}$ and $\Delta_D$ is the Dirichlet-Laplacian on $\rR^2\backslash\rK$.
When $s$ is an integer, such norm is equivalent to $\dot W^{s,2}$ norm. See \cite{MR2584618} for the detailed discussion of such space.
\end{remark}
\begin{remark}
Remark that $\rK$ is star-shaped, i.e. for any $x\in\partial\rK$ with outer normal vector $N_\rK$, there is 
$x\cdot N_\rK \geq 0$. 
By the idea in Farah \cite{MR3104079,MR3109476}, star-shaped condition can be changed to that so called illuminate condition.
Though, for the clarity of this paper, we will not discuss those situations.
\end{remark}

The rest of paper is organized as follows.
In the next section, we give a quick sketch of the corresponding local and global existence to \eqref{Eq:oe}.
In the third section, we give the proof of \Tm{Tm:ped} by using the energy method.
In the fourth section, we give the proof of \Tm{Tm:pd} by applying a kind of uniform Gagliardo-Nirenberg inequality.
In the last section, we establish a mixed time-space estimate and give the proof of \Tm{Tm:sr}.

In what follows, we use the notation $S_{T_1}^{T_2}:=[T_1,T_2]\times \rR^2\backslash\rK$ and $\phi^p:=|\phi|^{p-1}\phi$.
We omit the integral region or volume element while they do not cause confusion.
We also use the Einstein summation convention, 
as well as the convention that Greek indices $\mu, \nu, \cdots$ range from $0$ to $2$ while Latin indices $i, j, \cdots$ will run from $1$ to $2$.
Finally, without loss of generality, all the proof will only be considered in the future $\rR_+\times\rR^2\backslash\rK$. 

\section{Local and global existence}\label{Sn:Lage}
The local well posed of strong solution to \eqref{Eq:oe} in $H^2\times H^1$ follows easily from the energy inequality,
which proof is almost the same as that in $\rR^{1+2}$ (see, e.g., Evans \cite{MR2597943}).
The same approach can also  be applied in our case when $\varE_{1,0}$ is finite.

The proof for global existence is also fairly standard. 
For the completeness of this paper, we give a quick sketch.
Assume that the solution exists up to $T_*<\infty$.
We use the energy inequality in $S_0^{T_*}$ and get
\begin{Eq*}
\|\partial^2 \phi\|_{L_t^\infty L_x^2(S_0^{T_*})}\lesssim& \|\partial^2 \phi(0)\|_{L_x^2}+\int_{0}^{T_*}\kl\||\phi|^{p-1}|\partial \phi|\kr\|_{L_x^2}\d t\\
\lesssim& \|\partial^2 \phi(0)\|_{L_x^2}+T_*\|\phi\|_{L_t^\infty L_x^{2p}(S_0^{T_*})}^{p-1}\|\partial \phi\|_{L_t^\infty L_x^{2p}(S_0^{T_*})}.
\end{Eq*}
Here, using the {Sobolev} interpolation inequality, we know
\begin{Eq*}
\|\phi\|_{L_x^{2p}}\leq & \|\phi\|_{L_x^{p+1}}^{\frac{p+1}{2p}}\|\phi\|_{\dot H^{1}}^{\frac{p-1}{2p}}\leq  \|\phi\|_{L_x^{p+1}}^{\frac{p+1}{2p}}\|\partial\phi\|_{L_x^2}^{\frac{p-1}{2p}},\\
\|\partial\phi\|_{L_x^{2p}}\leq & \|\partial\phi\|_{L_x^{2}}^{\frac{1}{p}}\|\partial\phi\|_{\dot H^{1}}^{\frac{p-1}{p}}\leq  \|\partial\phi\|_{L_x^2}^{\frac{1}{p}}\|\partial^2\phi\|_{L_x^2}^{\frac{p-1}{p}},
\end{Eq*}
respectively. 
Noticing $\partial_t^2\phi=\Delta_x\phi-\phi^p$, 
we have 
\begin{Eq*}
\|\partial^2 \phi(0)\|_{L_x^2}\lesssim \|\partial_x\partial \phi(0)\|_{L_x^2}+\|\phi(0)\|_{L_x^{p+1}}^{\frac{p+1}{2}}\|\partial\phi(0)\|_{L_x^2}^{\frac{p-1}{2}}\leq C(\varE_{1,0}).
\end{Eq*}
Now, using the energy conservation law, we easily get that
\begin{Eq*}
\|\partial^2 \phi\|_{L_t^\infty L_x^2(S_0^{T_*})}\leq \frac{1}{2}\|\partial^2 \phi\|_{L_t^\infty L_x^2(S_0^{T_*})}+C(T_*,\varE_{1,0}).
\end{Eq*}
Thus, $\|\partial^2\phi\|_{L_x^2}$, and similarly $\|\partial^{\leq 2}\phi\|_{L_x^2}$, can not blow up on $t=T_*$.
This is in contradiction to the assumption that the solution exists up to $T_*$, which means that the solution must be global.

\section{Decay of the potential energy}%Dotpe
To begin with, let us give a sketch of the energy method.
Fixing a vector field $X$ and a corresponding function $\chi$, we define
\begin{Eq*}
T_{\mu\nu}:=&\partial_\mu \phi\partial_\nu \phi-\frac{1}{2}m_{\mu\nu}\kl(\partial^\sigma \phi\partial_\sigma \phi+\frac{2}{p+1}|\phi|^{p+1}\kr),\\
J_\mu:=&T_{\mu\nu} X^\nu-\frac{1}{2}\partial_\mu\chi\cdot|\phi|^2+\frac{1}{2}\chi\partial_\mu|\phi|^2,
\end{Eq*}
where $T_{\mu\nu}$ is called the energy momentum tensor and $J_\mu$ is called the current. 
Here we take $x^0=t$  and indices are raised and lowered with respect to the Minkowski metric $m=\diag(-1,1,1,1)$.
For any domain $\rD$, applying the {Stokes}' formula, we derive the energy identity
\begin{Eq*}
\int_{\partial \rD} J_{N_\rD}=\int_{\rD} \partial^\mu J_\mu
\end{Eq*}
where $N_\rD$ is its outer normal vector on $\partial\rD$ and $J_{N_\rD}:=J_\mu N_{\rD}^\mu$.
Here we choose $\rD=S_0^T$ so that
\begin{Eq*}
\int_{t=T} J_0+\int_{S_0^T} \partial^\mu J_\mu= \int_{t=0} J_0-\int_{[0,T]\times \partial\rK} J_{N}, \quad (N:=(0,N_\rK)~on~\rR\times \partial\rK).
\end{Eq*} 
We expect that $J_{N}\geq 0$ so that the quantities on the left hand side can be controlled by the initial data. 

Noticing $\phi=\partial_t\phi=0$ on $\rR\times \partial\rK$, we have
\begin{Eq*}
|\partial_x\phi|=|N\phi|,\qquad X\phi=(X^i N_i) N\phi.
\end{Eq*}
Thus, on $\rR\times \partial\rK$, we see
\begin{Eq*}
J_N=T_{\mu\nu} N^\mu X^\nu =N^iX^j\kl(\partial_i\phi\partial_j\phi-\frac{1}{2}m_{ij}|\partial_x\phi|^2\kr)=\frac{1}{2}(X^iN_i)|N\phi|^2.
\end{Eq*}
This means that to reach $J_N\geq 0$, we only need to ensure $N\cdot X\geq 0$, 
which then means that we need to find the suitable $X$ that works on $\rR^{1+2}$ and has expression
$X=X^t\partial_t+X^r\partial_r$ with $X^r\geq 0$, since that $\rK$ is star-shaped.

\begin{proof}[Proof of \Tm{Tm:ped}]
For $p\geq 5$, we set $X=(r^2+t^2+1)\partial_t+2tr\partial_r$ with $\chi=t$. 
Obviously, we have $X^r\geq 0$. 
So, similar to that in $\rR^{1+2}$, we know
\begin{Eq*}
&\tilde E(T)+\frac{p-5}{p+1}\int_{S_0^T}t|\phi|^{p+1}\leq \tilde E(0)\lesssim \varE_{0,2},\\
&\tilde E(t):=\int_{\rR^2\backslash\rK}\frac{1}{2}\kl(|t\partial_t\phi+r\partial_r+\phi|^2+\sum_{\mu<\nu}|x_\mu\partial_\nu\phi-x_\nu\partial_\mu\phi|^2+|\partial\phi|^2\kr)+\frac{t^2+r^2+1}{p+1}|\phi|^{p+1}.
\end{Eq*}
The process of its proof is standard. 
We recommend interested readers to, e.g., Glassey and Pecher \cite{MR667924} and Alinhac \cite{MR2524198} to get the calculation details in $\rR^{1+2}$, and modify them to $\rR\times\rR^2\backslash\rK$.
As a result, we get \eqref{Eq:ped2}.
 
For $1<p<5$, we first consider the $\rR^{1+2}$ case
and introduce the two vector fields constructed in \cite{MR4395159},
\begin{Eq*}
\tilde X_1:=&(x_2^2+(t-x_1)^2+1)\partial_t+(x_2^2-(t-x_1)^2)\partial_1+2(t-x_1)x_2\partial_2,\\
\tilde X_2:=&(t-x_1+1)^{\frac{p-1}{2}}(\partial_t-\partial_1)+(t-x_1+1)^{\frac{p-5}{2}}x_2^2(\partial_t+\partial_1)+2(t-x_1+1)^{\frac{p-3}{2}}x_2\partial_2.
\end{Eq*}
For $\tilde X_1$ with $\tilde \chi_1:=t-x_1$, there are
\begin{Eq*}
\tilde J_{1;\partial_t}=&\tilde J_{1;\partial_t}^M+\frac{1}{2}\partial_1((t-x_1)|\phi|^2)-\frac{1}{2}\partial_2(x_2|\phi|^2),\\
\tilde J_{1;\partial_t}^M:=&\frac{1}{2}\bigg(|x_2(\partial_t\phi+\partial_1\phi)+(t-x_1)\partial_2\phi|^2+|(t-x_1)(\partial_t\phi-\partial_1\phi)+x_2\partial_2\phi+\phi|^2\\
&\phantom{\frac{1}{2}\bigg(}+|\partial_t\phi|^2+|\partial_x\phi|^2+\frac{2}{p+1}\kl((t-x_1)^2+x_2^2+1\kr)|\phi|^{p+1}\bigg),\\
\kl.\tilde J_{1;\partial_t+\partial_1}\kr|_{t=x_1}=&\kl(x_2^2+\frac{1}{2}\kr)|\partial_t\phi+\partial_1\phi|^2+\frac{1}{2}|\partial_2\phi|^2+\frac{1}{p+1}|\phi|^{p+1},\\
\partial^{\mu}\tilde J_{1;\mu}=&\frac{5-p}{p+1}(x_1-t)|\phi|^{p+1}.
\end{Eq*}
Here we use its corresponding energy method on $\{0\leq t\leq \min\{T,x_1\}\}$,
where $\tilde J_{1;\partial_t}^M$, $\kl.\tilde J_{1;\partial_t+\partial_1}\kr|_{t=x_1}$ and $\partial^{\mu}\tilde J_{1;\mu}$ are nonnegative, 
and the rest terms in $\tilde J_{1;\partial_t}$ vanish in the integral.
As a consequence, we especially have
\begin{Eq*}
\int_{0\leq t=x_1\leq T}\tilde J_{1;\partial_t+\partial_1}\leq \int_{0\leq x_1}\kl.\tilde J_{1;\partial_t}^M\kr|_{t=0}\lesssim\varE_{0,2}.
\end{Eq*}
Meanwhile, for $\tilde X_2$ with $\tilde \chi_2:=(t-x_1+1)^{(p-1)/2}$, there are
\begin{Eq*}
\tilde J_{2;\partial_t}=&\tilde J_{2;\partial_t}^M+\frac{1}{2}\partial_1\kl((t-x_1+1)^{\frac{p-3}{2}}|\phi|^2\kr)-\frac{1}{2}\partial_2\kl(x_2(t-x_1+1)^{\frac{p-5}{2}}|\phi|^2\kr),\\
\tilde J_{2;\partial_t}^M:=&\frac{1}{2}\bigg(|x_2(\partial_t\phi+\partial_1\phi)+(t-x_1+1)\partial_2\phi|^2+|(t-x_1+1)(\partial_t\phi-\partial_1\phi)+x_2\partial_2\phi+\phi|^2\\
&\phantom{\frac{1}{2}\bigg(}+\frac{2}{p+1}\kl((t-x_1+1)^2+x_2^2\kr)|\phi|^{p+1}\bigg)(t-x_1+1)^{\frac{p-5}{2}},\\
\kl.\tilde J_{2;\partial_t+\partial_1}\kr|_{t=x_1}=&\tilde J_{2;\partial_t+\partial_1}^M+\frac{1}{2}(\partial_t+\partial_1)|\phi|^2\\
\tilde J_{2;\partial_t+\partial_1}^M:=&|x_2(\partial_t\phi+\partial_1\phi)+\partial_2\phi|^2+\frac{2}{p+1}|\phi|^{p+1},\\
\partial^{\mu}\tilde J_{2;\mu}=&\frac{5-p}{2}(t-x_1+1)^{\frac{p-7}{2}}|x_2(\partial_t\phi+\partial_1\phi)+(t-x_1+1)\partial_2\phi|^{2}.
\end{Eq*}
Here we use its corresponding energy method on $\{\max\{0,x_1\}\leq t\leq T\}$,
where $\tilde J_{2;\partial_t}^M$, $\tilde J_{2;\partial_t+\partial_1}^M$ and $\partial^{\mu}\tilde J_{2;\mu}$  are nonnegative, 
and the rest terms in $\tilde J_{2;\partial_t}$ and $\kl.\tilde J_{2;\partial_t+\partial_1}\kr|_{t=x_1}$ all vanish in the integral.
As a consequence, we especially have
\begin{Eq*}
\int_{x_1\leq T}\kl.\tilde J_{2;\partial_t}^M\kr|_{t=T}\leq \int_{x_1\leq 0}\kl.\tilde J_{2;\partial_t}^M\kr|_{t=0}+\int_{0\leq t=x_1\leq T}\tilde J_{2;\partial_t+\partial_1}^M.
\end{Eq*}
Finally, noticing $\tilde J_{2;\partial_t+\partial_1}^M\leq 4 \tilde J_{1;\partial_t+\partial_1}$, for equation on $\rR^{1+2}$, we are able to get
\begin{Eq}\label{Eq:FroR}%Final result on R^{1+2}
\int_{\rR^{1+2}}{\bf 1}_{x_1\leq 0}\cdot(t-x_1+1)^{\frac{p-1}{2}}|\phi|^{p+1}\d x\lesssim \varE_{0,2}.
\end{Eq}

Now, mixing these two vector fields, we define 
\begin{Eq*}
\tilde X:={\bf 1}_{t\leq x_1}\cdot 4\tilde X_1+{\bf 1}_{t> x_1}\cdot \tilde X_2,\qquad \tilde \chi:={\bf 1}_{t\leq x_1}\cdot 4\tilde \chi_1+{\bf 1}_{t> x_1}\cdot \tilde \chi_2.
\end{Eq*}
The only obstacle to using it in our exterior problem is that such vector field is not spherically symmetric.
For this reason, we construct the spherical integral of $\tilde X$.
To be more specific, we use polar coordinates and write
\begin{Eq*}
\tilde X(t,r,\theta)=\tilde X^t(t,r,\theta)\partial_t+\tilde X^r(t,r,\theta)\partial_r+\tilde X^\theta(t,r,\theta)\partial_\theta.
\end{Eq*}
Then we define
\begin{Eq*}
\chi(t,r,\theta):=&\int_0^{2\pi} \tilde\chi(t+R,r,\theta+\theta')\d\theta',\\
X(t,r,\theta):=&\int_0^{2\pi} \tilde X(t+R,r,\theta+\theta')\d\theta'\\
=&\kl(\int_0^{2\pi} \tilde X^t(t+R,r,\theta')\d\theta'\kr) \partial_t+\kl(\int_0^{2\pi} \tilde X^r(t+R,r,\theta')\d\theta'\kr) \partial_r.
\end{Eq*}
Here we mention that $\int_0^{2\pi} \tilde X^\theta(\theta')\d\theta'=0$ because $\tilde X^\theta$ is an odd function, 
which in turn is due to the fact that $\tilde X$ is symmetric with respect to $x_2$. 
Using its corresponding energy method on $[0,t]\times \rR^2$
and investigating by \eqref{Eq:FroR}, we are able to get that
\begin{Eq*}
\varE_{0,2}\gtrsim&\int_0^{2\pi}\int_{\rR^2} {\bf 1}_{\cos(\theta+\theta')\leq 0}(t-r\cos(\theta+\theta')+R+1)^{\frac{p-1}{2}}|\phi|^{p+1}(t,r\cos\theta,r\sin\theta)\d x\d\theta'\\
\gtrsim& \int_{\rR^2}\int_{-\theta+5\pi/6}^{-\theta+7\pi/6}(t-r\cos(\theta+\theta')+R+1)^{\frac{p-1}{2}}|\phi|^{p+1}\d\theta'\d x\\
\gtrsim&\int_{\rR^2} \kl(t+\frac{r}{2}+R+1\kr)^{\frac{p-1}{2}}|\phi|^{p+1}\d x
\gtrsim\int_{\rR^2} \kl(t+r+1\kr)^{\frac{p-1}{2}}|\phi|^{p+1}\d x.
\end{Eq*}

On the other hand, on $\rR_+\times\partial\rK$, we always have $t+R\geq r$ and thus
\begin{Eq*}
X^r=&\int_0^{2\pi} \tilde X^r(t+R,r,\theta')\d\theta'=\int_0^{2\pi} \tilde X_2^r(t+R,r,\theta)\d\theta\\
=&\int_0^{2\pi} \kl(-(t-r\cos\theta+R+1)^{\frac{p-1}{2}}+(t-r\cos\theta+R+1)^{\frac{p-5}{2}}(r\sin\theta)^2\kr)\cos\theta\d\theta\\
&+\int_0^{2\pi}2(t-r\cos\theta+R+1)^{\frac{p-3}{2}}r\sin\theta\sin\theta\d\theta\\
=&2\int_{0}^{\frac{\pi}{2}} \kl(-(t-r\cos\theta+R+1)^{\frac{p-1}{2}}+(t+r\cos\theta+R+1)^{\frac{p-1}{2}}\kr)\cos\theta\d\theta\\
&+2\int_{0}^{\frac{\pi}{2}} \kl((t-r\cos\theta+R+1)^{\frac{p-5}{2}}-(t+r\cos\theta+R+1)^{\frac{p-5}{2}}\kr)r^2\cos\theta\sin^2\theta\d\theta\\
&+4\int_{0}^{\frac{\pi}{2}} \kl((t-r\cos\theta+R+1)^{\frac{p-3}{2}}+(t+r\cos\theta+R+1)^{\frac{p-3}{2}}\kr)r\sin^2\theta\d\theta\\
\geq&0.
\end{Eq*}
Through the discussion at the beginning of this section,
we finish the proof of the \Tm{Tm:ped}.
\end{proof}

\section{Pointwise decay of the solution}

To begin with, we introduce the necessary bound for each kind of energy.
\begin{proposition}\label{Pn:ked}%kinds energy decay
Assume $\phi$ is the solution of \eqref{Eq:oe}. We have
\begin{Eq*}
(1+t)^2&\|\spartial \phi\|_{L^2(r>R)}^2+\|(1+|t-r|)\partial \phi\|_{L_x^2}^2+\|\phi\|_{L_x^2}^2\lesssim \begin{cases}
(1+t)^{\frac{5-p}{2}}, &p<5,\\
1, &p\geq 5,
\end{cases}
\end{Eq*}
and $\|\phi\|_{H^2}\lesssim (1+t)^2$, where the constants of inequalities depends on $\varE_{0,2}$ and $\varE_{1,0}$.
\end{proposition}
The proof of \Pn{Pn:ked} is identical to the proof of Proposition 4.2 in \cite{MR4395159}. 
The star-shape of $\rK$ joint with the Dirichlet boundary condition ensure that all the process can be adapted here. 
We only mention that $\spartial\phi$ is formally defined in whole $\rR^2\backslash\rK$ but only be useful away from $\rK$,
and that the decay of energy is the same as that of homogeneous solution when $p\geq5$.

On $\rR^{1+2}$, such estimate joint with the suitable version of {Sobolev} inequality gives the pointwise decay of solution.
See \cite{MR4395159} for the detailed discussion.
Thus, on $\rR\times\rR^2\backslash\rK$, the main difficulty will arise near $\partial\rK$. 
\begin{proof}[Proof of \Tm{Tm:pd}]
Without loss of generality, we consider $t+r\geq 10R$,
and divide this region into three sub regions: 
$\{r\leq t/2\}$, $\{t/2\leq r\leq 3t/2\}$ and $\{3t/2\leq r\}$.
The proof of estimates in the last two regions is almost the same as that in \cite{MR4395159}
since these regions are far away from $\rR\times \partial\rK$.
So, we only need to deal with the first region.

For any fixed $t\geq 5R$, we define 
\begin{Eq*}
\tilde \phi(x):=\phi(t,tx)~for~x\in B_{1/2}\backslash t^{-1}\rK,
\end{Eq*}
where $t^{-1}\rK=\{t^{-1}x:x\in\rK\}$.
In view of \Pn{Pn:ked}, we conclude
\begin{Eq*}
\int_{B_{1/2}\backslash (t^{-1}\rK)}|\tilde\phi|^2\leq& t^{-2}\int_{r<\frac{1}{2}t}|\phi|^2\lesssim \begin{cases}
(1+t)^{\frac{1-p}{2}},&p<5,\\
(1+t)^{-2},&p\geq 5,
\end{cases}\\
\int_{B_{1/2}\backslash (t^{-1}\rK)}|\partial_x\tilde\phi|^2\leq& \int_{r<\frac{1}{2}t}|\partial_x\phi|^2\lesssim t^{-2}\int_{r<\frac{1}{2}t}|(t-r)\partial_x\phi|^2\lesssim 
\begin{cases}
(1+t)^{\frac{1-p}{2}},&p<5,\\
(1+t)^{-2},&p\geq5,
\end{cases}\\
\int_{B_{1/2}\backslash (t^{-1}\rK)}|\partial_x^2\tilde\phi|^2\leq& t^2\int_{r<\frac{1}{2}t}|\partial_x^2\phi|^2\lesssim (1+t)^6
\end{Eq*}
separately. 
Next, we are hoping to use the inequality of embedding.
To do so, we need to find a uniform constant in the inequality 
though the region $B_{1/2}\backslash (t^{-1}\rK)$ depends on $t$.

For this purpose, we first consider $B_{1/2}\backslash (5R)^{-1}\rK$.
Noticing $B_{1/2}$ is convex and $(5R)^{-1}\rK$ is star-shaped with smooth boundary, we can easily construct a triangle $A$ such that for any $x\in B_{1/2}\backslash (5R)^{-1}\rK$, there exists an $A_x$ satisfying
\begin{Eq*}
x\in A_x,\qquad A_x\cong A,\qquad A_x\subset B_{1/2}\backslash (5R)^{-1}\rK.
\end{Eq*} 
Moreover, for such $A$ and any $t>5R$, it is easy to find that for any $x\in B_{1/2}\backslash t^{-1}\rK$, 
there also exists an $A_x$ satisfying
\begin{Eq*}
x\in A_x,\qquad A_x\cong A,\qquad A_x\subset B_{1/2}\backslash t^{-1}\rK.
\end{Eq*}
See the figure below for a simple geometric intuition of this statement.
\begin{figure}[H]
\centering
\includegraphics[width=0.7\textwidth]{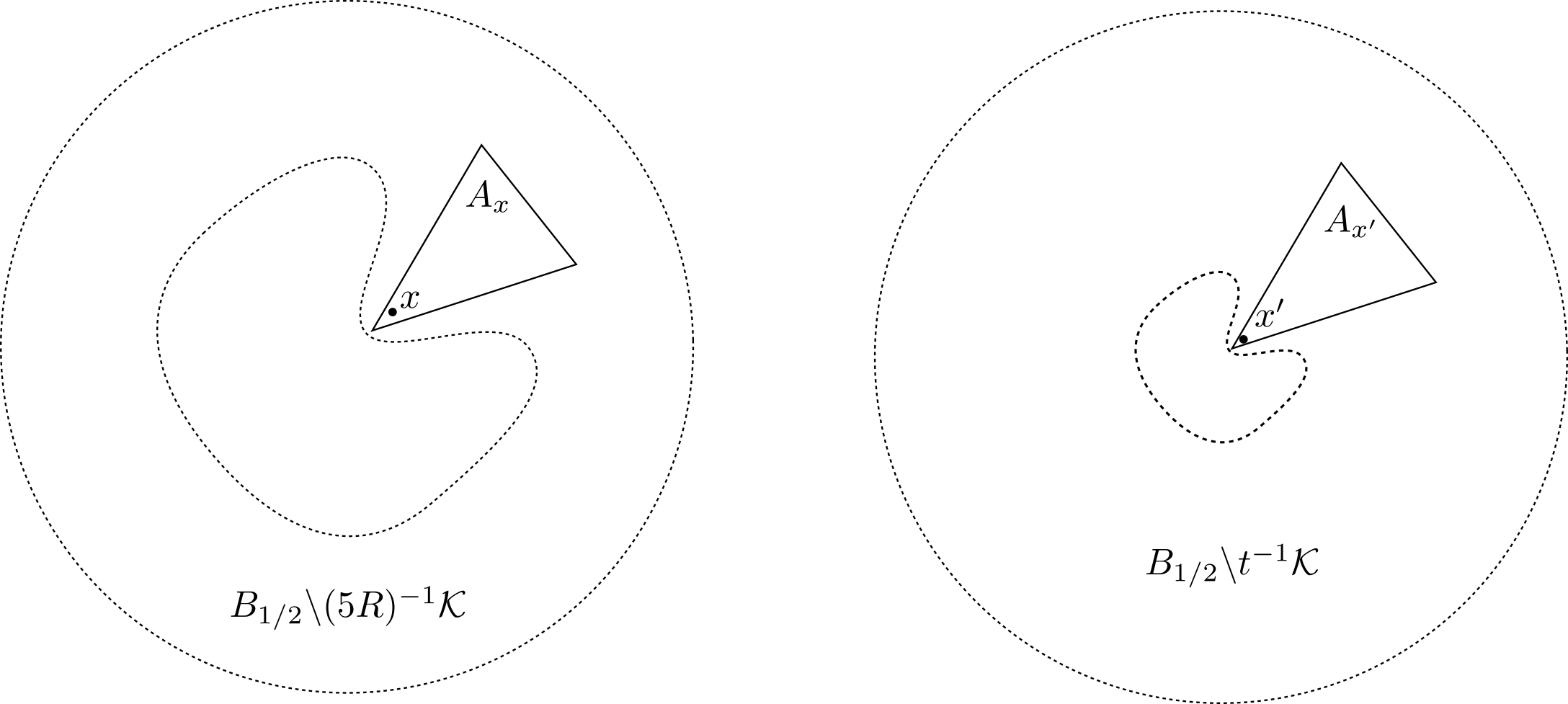}
\end{figure}
Now, for any fixed $x_0\in B_{1/2}\backslash t^{-1}\rK$, using the {Gagliardo-Nirenberg} inequality, we see
\begin{Eq*}
\|\tilde\phi\|_{L_x^\infty(A_{x_0})}\lesssim& \|\tilde\phi\|_{L_x^2(A_{x_0})}^{\delta}\|\tilde\phi\|_{\dot W^{1,\frac{2-2\delta}{1-2\delta}}(A_{x_0})}^{1-\delta}\\
\lesssim& \|\tilde\phi\|_{L_x^2(A_{x_0})}^{\delta}\|\tilde\phi\|_{\dot W^{1,2}(A_{x_0})}^{1-2\delta}\|\tilde\phi\|_{\dot W^{2,2}(A_{x_0})}^{\delta}\\
\lesssim&(1+t)^{\frac{1-\min\{p,5\}}{4}+\frac{\min\{p,5\}+23}{4}\delta},
\end{Eq*}
where the constant only depends on $A$ but not on $t$. Now, choosing $\delta$ sufficiently small with respect to $\varepsilon$, we see
\begin{Eq*}
\|\phi\|_{L_x^\infty(r<t/2)}=&\|\tilde\phi\|_{L_x^\infty(B_{1/2}\backslash t^{-1}\rK)}
=\sup_{x_0\in B_{1/2}\backslash t^{-1}\rK} \|\tilde\phi\|_{L_x^\infty(A_{x_0})}\\
\lesssim&(1+t)^{\frac{1-\min\{p,5\}}{4}+\varepsilon},
\end{Eq*}
which finishes the proof of \Tm{Tm:pd}.
\end{proof}

\section{Scattering of the solution}
Firstly, we prove the scattering result in energy space.
For writing convenience, for any $s$, we define
\begin{Eq*}
{\bf f}:=(f,\partial_t f),\qquad
\|{\bf f}\|_{\dot\rH^s}:=\|f\|_{\dot H^s}+\|\partial_tf\|_{\dot H^{s-1}}.
\end{Eq*}
\begin{proof}[Proof of \eqref{Eq:sr1}]

By energy inequality, for any $T_2>T_1$, we can calculate that
\begin{Eq*}
\|\boldsymbol\phi(T_2)-L(T_2-T_1)\boldsymbol\phi(T_1)\|_{\dot\rH^1}\lesssim \|\phi^p\|_{L_t^1L_x^2(S_{T_1}^{T_2})}\lesssim \|\phi\|_{L_t^pL_x^{2p}(S_{T_1}^{T_2})}^p.
\end{Eq*}
Using the interpolation inequality, H\"older inequality, \Tm{Tm:ped} and \Tm{Tm:pd}, for $2\sqrt{5}-1<p<5$ we see
\begin{Eq*}
\|\phi\|_{L_x^{2p}}
\lesssim& \|\phi\|_{L_x^{p+1}}^{\frac{p+1}{2p}}\|\phi\|_{L_x^\infty}^{\frac{p-1}{2p}}\\
\lesssim& \kl\|(1+t+r)^{\frac{p-1}{2}}|\phi|^{p+1}\kr\|_{L_x^{1}}^{\frac{1}{2p}}\kl\|(1+t+r)^{-\frac{p-1}{2}}\kr\|_{L_x^\infty}^{\frac{1}{2p}}\|\phi\|_{L_x^\infty}^{\frac{p-1}{2p}}\\
\lesssim&(1+t)^{\frac{1}{2p}(-\frac{p-1}{2})+\frac{p-1}{2p}(-\frac{p-1}{8}+\varepsilon)}\\
\lesssim& (1+t)^{-\frac{p^2+2p-3}{16p}+\frac{p-1}{2p}\varepsilon}.
\end{Eq*}
Then, noticing $-\frac{p^2+2p-3}{16}+\frac{p-1}{2}\varepsilon<-1$ while $p>2\sqrt5-1$ and $\varepsilon$  small enough, we have
\begin{Eq*}
 \|\phi\|_{L_t^pL_x^{2p}(S_{T_1}^{T_2})}^p\lesssim& \int_{T_1}^{T_2}(1+t)^{-\frac{p^2+2p-3}{16}+\frac{p-1}{2}\varepsilon}d t\lesssim (1+T_1)^{-\frac{p^2+2p-19}{16}+\frac{p-1}{2}\varepsilon}.
\end{Eq*}

Now, we easily see that
\begin{Eq*}
\sup_{T\geq 0}\|L(-T)\boldsymbol\phi(T)\|_{\dot\rH^{1}}\leq& \|\boldsymbol\phi(0)\|_{\dot\rH^{1}}+\sup_{T\geq 0}\|\boldsymbol\phi(T)-L(T)\boldsymbol\phi(0)\|_{\dot\rH^{1}}<\infty,\\
\lim_{T_1,T_2\rightarrow \infty}\|L(-T_2)\boldsymbol\phi(T_2)-L(-T_1)\boldsymbol\phi(T_1)\|_{\dot\rH^{1}}=&\lim_{T_1,T_2\rightarrow \infty}\|\boldsymbol\phi(T_2)-L(T_2-T_1)\boldsymbol\phi(T_1)\|_{\dot\rH^{1}}= 0.
\end{Eq*}
Defining $(\phi_{0}^+(x),\phi_1^+(x))=\lim_{T\rightarrow\infty}L(-T)\boldsymbol\phi(T)$ in $\rH^1$, we get \eqref{Eq:sr1} for these cases of $p$. 

Next, for $p\geq 5$, we similarly get
\begin{Eq*}
\|\phi\|_{L_x^{2p}}\lesssim& \|\phi\|_{L_x^{p+1}}^{\frac{p+1}{2p}}\|\phi\|_{L_x^\infty}^{\frac{p-1}{2p}}\lesssim (1+t)^{\frac{p+1}{2p}(-\frac{2}{p+1})+\frac{p-1}{2p}(-\frac{1}{2}+\varepsilon)}\lesssim (1+t)^{-\frac{p+3}{4p}+\frac{p-1}{2p}\varepsilon},\\
 \|\phi\|_{L_t^pL_x^{2p}(S_{T_1}^{T_2})}^p\lesssim& \int_{T_1}^{T_2}(1+t)^{-\frac{p+3}{4}+\frac{p-1}{2}\varepsilon}d t\lesssim (1+T_1)^{-\frac{p-1}{4}+\frac{p-1}{2}\varepsilon}.
\end{Eq*}
Through the same process, we finish the proof.
\end{proof}

\begin{proof}[Proof of \eqref{Eq:sr2}]
Similar to before, we only need to show
$\|\boldsymbol\phi(T_2)-L(T_2-T_1)\boldsymbol\phi(T_1)\|_{\dot\rH^{s}}$ converges to $0$ as $T_1$ tends to infinity and uniformly 
for $T_2>T_1$.
Also, by interpolation, we only need to consider the critical case $s=s_p=\frac{p-3}{p-1}$.

We first adopt the $\rH^s$ energy inequality and get
\begin{Eq*}
\|\boldsymbol\phi(T_2)-L(T_2-T_1)\boldsymbol\phi(T_1)\|_{\dot\rH^{s_p}}
\lesssim& \|\phi^p\|_{L_t^1\dot H_x^{s_p-1}(S_{T_1}^{T_2})}
\lesssim\|\phi\|_{L_t^p L_x^{\frac{2p(p-1)}{p+1}}(S_{T_1}^{T_2})}^p.
\end{Eq*}

When $1+2\sqrt{2}<p\leq 2+\sqrt{5}$, we have $2p(p-1)/(p+1)\leq p+1$.
Then, using H\"older inequality and \Tm{Tm:ped}, we see
\begin{Eq*}
\|\phi\|_{L_x^{\frac{2p(p-1)}{p+1}}}\lesssim&\kl\|(1+t+r)^{\frac{p-1}{2}}|\phi|^{p+1}\kr\|_{L_x^{1}}^{\frac{1}{p+1}}\kl\|(1+t+r)^{-\frac{p-1}{2}}\kr\|_{L_x^\frac{2p(p-1)}{-p^2+4p+1}}^{\frac{1}{p+1}}\\
\lesssim& (1+t)^{-\frac{p^3-7p-2}{2p(p^2-1)}}.
\end{Eq*}
Noticing $-\frac{p^3-7p-2}{2(p^2-1)}<-1$ while $p>2\sqrt2+1$, we have
\begin{Eq*}
\|\phi\|_{L_t^pL_x^{\frac{2p(p-1)}{p+1}}(S_{T_1}^{T_2})}^p\lesssim& \int_{T_1}^{T_2} (1+t)^{-\frac{p^3-7p-2}{2(p^2-1)}}\d t\lesssim (1+T_1)^{-\frac{p(p^2-2p-7)}{2(p^2-1)}}
\end{Eq*}
and finish the proof for these cases of $p$.

When $2+\sqrt{5}<p< 5$, where $2p(p-1)/(p+1)> p+1$,
similar to the proof of \eqref{Eq:sr1},
we see
\begin{Eq*}
\|\phi\|_{L_x^{\frac{2p(p-1)}{p+1}}}
%\lesssim& \|\phi\|_{L_x^{p+1}}^{\frac{(p+1)^2}{2p(p-1)}}\|\phi\|_{L_x^\infty}^{\frac{p^2-4p-1}{2p(p-1)}}\\
%\lesssim&\||\phi|^{p+1}\|_{L_x^1}^{\frac{p+1}{2p(p-1)}}\|\phi\|_{L_x^\infty}^{\frac{p^2-4p-1}{2p(p-1)}}\\
\lesssim&\kl\|(1+t+r)^{\frac{p-1}{2}}|\phi|^{p+1}\kr\|_{L_x^1}^{\frac{p+1}{2p(p-1)}}\kl\|\kl(1+t+r\kr)^{-\frac{p-1}{2}}\kr\|_{L_x^\infty}^{\frac{p+1}{2p(p-1)}}\|\phi\|_{L_x^\infty}^{\frac{p^2-4p-1}{2p(p-1)}}\\
%\lesssim&(1+t)^{\frac{p+1}{2p(p-1)}(-\frac{p-1}{2})+\frac{p^2-4p-1}{2p(p-1)}(-\frac{p-1}{8}+\varepsilon)}\\
\lesssim&(1+t)^{-\frac{p^2+3}{16p}+\frac{p^2-4p-1}{2p(p-1)}\varepsilon}.
\end{Eq*}
Noticing $-\frac{p^2+3}{16}+\frac{p^2-4p-1}{2(p-1)}\varepsilon<-1$ while 
$2+\sqrt{5}<p< 5$ and $\varepsilon$ small enough, we have
\begin{Eq*}
\|\phi\|_{L_t^pL_x^{\frac{2p(p-1)}{p+1}}(S_{T_1}^{T_2})}^p\lesssim& \int_{T_1}^{T_2} (1+t)^{-\frac{p^2+3}{16}+\frac{p^2-4p-1}{2(p-1)}\varepsilon}\d t\\
\lesssim& (1+T_1)^{-\frac{p^2-13}{16}+\frac{p^2-4p-1}{2(p-1)}\varepsilon}
\end{Eq*}
and finish the proof for these cases of $p$.

Finally, for $ p\geq 5$, similar to above, we see
\begin{Eq*}
\|\phi\|_{L_x^{\frac{2p(p-1)}{p+1}}}
\lesssim&\kl\|(1+t+r)^{2}|\phi|^{p+1}\kr\|_{L_x^1}^{\frac{p+1}{2p(p-1)}}\kl\|\kl(1+t+r\kr)^{-2}\kr\|_{L_x^\infty}^{\frac{p+1}{2p(p-1)}}\|\phi\|_{L_x^\infty}^{\frac{p^2-4p-1}{2p(p-1)}}\\
%\lesssim&(1+t)^{\frac{p+1}{2p(p-1)}(-2)+\frac{p^2-4p-1}{2p(p-1)}(-\frac{1}{2}+\varepsilon)}\\
\lesssim&(1+t)^{-\frac{p^2+3}{4p(p-1)}+\frac{p^2-4p-1}{2p(p-1)}\varepsilon},\\
\|\phi\|_{L_t^pL_x^{\frac{2p(p-1)}{p+1}}(S_{T_1}^{T_2})}^p\lesssim& \int_{T_1}^{T_2} (1+t)^{-\frac{p^2+3}{4(p-1)}+\frac{p^2-4p-1}{2(p-1)}\varepsilon}\d t\\
\lesssim& (1+T_1)^{-\frac{p^2-4p+7}{4(p-1)}+\frac{p^2-4p-1}{2(p-1)}\varepsilon},
\end{Eq*}
where $-\frac{p^2-4p+7}{4(p-1)}+\frac{p^2-4p-1}{2(p-1)}\varepsilon< 0$ while $p\geq 5$ and $\varepsilon$ small enough.
This finishes the proof of \Tm{Tm:sr}.
\end{proof}

\subsection*{Acknowledgment}
The authors would like to thank the anonymous referee for careful reading and valuable comments.

\end{document}